\documentclass[conference]{IEEEtran}
\IEEEoverridecommandlockouts
\usepackage{spconf,amsmath,graphicx}
\usepackage{amsfonts}
\usepackage{wrapfig}
\usepackage{amsthm}
\usepackage{amssymb}
\usepackage{float}
\usepackage{xcolor}
\usepackage{algorithm}
\usepackage{algpseudocode}
\ifCLASSINFOpdf
\else
\fi

\newtheorem{theorem}{Theorem}
\newtheorem{lemma}{Lemma}
\newtheorem{assumption}{Assumption}

\newcommand*{\QEDB}{\hfill\ensuremath{\square}}%
\title{Dynamic Power Allocation for Smart Grids via ADMM}
\name{Marie Maros, Joakim Jald\'{e}n\thanks{This project has received funding from the European Research Council (ERC) under the European Union's Horizon 2020 research and innovation program (Grant Agreement no. 742648)}}
\address{Department of Information Science and Engineering \\
Royal Institute of Technology (KTH) \\
SE-100 44 Stockholm, Sweden}
\begin{document}

\maketitle
\begin{abstract}
Electric power distribution systems encounter fluctuations in supply due to renewable sources with high variability in generation capacity.
It is therefore necessary to provide algorithms that are capable of dynamically finding approximate solutions. We propose two semi-distributed algorithms based on ADMM and discuss their advantages and disadvantages. One of the algorithms computes a feasible approximate of the optimal power allocation at each time instance. We require coordination between the nodes to guarantee feasibility of each of the iterates. We bound the distance from the approximate solutions to the optimal solution as a function of the variation in optimal power allocation, and  we verify our results via experiments.
\end{abstract}
\vspace{-5pt}
\section{Introduction}
\label{sec:intro}
The introduction of renewable sources in energy markets poses new challenges that affect the power allocation policies of distribution systems \cite{challenge1}, \cite{challenge2}, \cite{matrix}. Decisions concerning the operation of the system including the dispatch of generation, must be done quickly due to the high variability in generation capacity of renewable sources. Distributed solutions are of interest due to their scalability with the number of energy sources and loads. However, distributed solutions that tackle the dispatch of generation typically deal with static scenarios where the production costs (or utilities), production capacity and load do not change over time. These solutions are typically iterative \cite{EconDisp} and the power balance constraint is typically not fulfilled until convergence is reached. In the time varying case this implies that a feasible solution may never be reached. Therefore, fast converging algorithms are desired so as to minimize the impact of problem variability.

We consider, like \cite{tracking}, a resource allocation problem with $N$ users in which the system operator's goal is to maximize the social welfare subject to a given system resource. Unlike \cite{tracking} we admit in our analysis the presence of box constraints, which, in the Economic Dispatch case, may represent constraints on the generation capacities. The utilities, constraints and system resources are considered to vary at the same time scale as the algorithm, i.e., we are only allowed one iteration before the problem changes. This captures the supply variability of renewable sources while also allowing for changes in generation cost and constraints. Further, we allow for smart consumers that have utilities dependent on their own time dependent demand. In contrast to most of the current literature, we do require the presence of a master node that supplies a limited amount of coordination. However, this ensures that each of the obtained solutions is feasible at each iteration. While this will increase the communication over head per time instance, only one iteration will be performed per time step.  In short, we propose and provide guarantees for a scheme in which a semi-decentralized power allocation algorithm, based on ADMM, is implemented in real time. Our main contribution guarantees that the iterates will remain at a bounded neighborhood of the optimal point and that the size of this neighborhood vanishes if the problem stops changing.
\vspace{-6pt}
\section{Problem Formulation and Algorithms}
\vspace{-2pt}
\label{sec:pf}
We consider a power distribution system with $R$ primary suppliers and $N$ users. In particular, we consider that the power injected by the primary suppliers is given, at time $k,$ by $\mathbf{P}^{[k]} \in \mathbb{R}^{R}.$ Further, the users may locally produce power and sell it to the primary suppliers to be used by other users in the system. The goal is to maximize the system's aggregate utility while taking into account $\mathbf{P}^{[k]},$ and the maximum and minimum power, $\overline{\mathbf{p}_i}^{[k]}$ and $\underline{\mathbf{p}_i}^{[k]},$ a user can produce and or consume at a specific time. The system's aggregate utility is defined considering the users' cost to produce electricity $C_i(\mathbf{p}_i,k)$ and the utility they obtain by consuming electricity $U_i(\mathbf{p}_i,k).$ This can be formulated mathematically as
\vspace{-2pt}
\begin{subequations}
\label{eq:original}
\begin{align}
& \underset{\{\mathbf{p}_{i}\}_{i=1}^N}{\min} & \sum_{i=1}^{M} C_i(\mathbf{p}_i,k) - \sum_{i=M+1}^{N}U_i(\mathbf{p}_i,k) \\
&\text{s.t.} &\underline{\mathbf{p}}_i^{[k]} \leq \mathbf{p}_i \leq \overline{\mathbf{p}}_i^{[k]}, \, i =1,\hdots,N \\
& & \sum_{i=1}^N \mathbf{p}_i = \mathbf{P}^{[k]}, \label{eq:powerbalance}
\end{align}
\end{subequations}
where all the quantities indexed by $[k]$ are allowed to vary in time.
The classic Economic Dispatch problem \cite{EconDisp}, where one intends to minimize the aggregate cost of producing electricity while meeting a demand constraint, can also be formulated in the form of \eqref{eq:original}. 

ADMM can be applied to solve \eqref{eq:original} in at least two different ways.   These are based on the following two reformulations of \eqref{eq:original}.
\begin{subequations}
\label{eq:partial}
\begin{align}
& \underset{\{\mathbf{p}_i\}_{i=1}^N,\, \{\mathbf{q}_i\}_{i=1}^N}{\min} & \sum_{i=1}^{M} C_i(\mathbf{p}_i,k) -  \sum_{i=M+1}^N U_i(\mathbf{q}_i,k) \\
& \text{s.t.} & \underline{\mathbf{p}}_i^{[k]} \leq \mathbf{p}_i \leq \overline{\mathbf{p}}_i^{[k]}, \, i =1,\hdots,N \label{eq:diff1}\\
& \quad  \sum_{i=1}^N \mathbf{q}_i = \mathbf{P}^{[k]} &  \mathbf{q}_i = \mathbf{p}_i,\, i =1,\hdots,N. \label{eq:consistency_partial}
\end{align}
\end{subequations}
\begin{subequations}
\label{eq:total}
\begin{align}
& \underset{\{\mathbf{p}_i\}_{i=1}^N,\, \{\mathbf{q}_i\}_{i=1}^N}{\min} & \sum_{i=1}^M C_i(\mathbf{p}_i,k) - \sum_{i=M+1}^N U_i(\mathbf{q}_i,k) \\
&  \text{s.t.} & \underline{\mathbf{p}}_i^{[k]} \leq \mathbf{q}_i \leq \overline{\mathbf{p}}_i^{[k]},\, i = 1,\hdots,N \label{eq:diff2}\\
& \quad \sum_{i=1}^N \mathbf{q}_i = \mathbf{P}^{[k]} 
&  \mathbf{q}_i = \mathbf{p}_i,\, i=1,\hdots,N.\label{eq:consistency}
\end{align}
\end{subequations}
The reformulation we use may impact the algorithm's convergence rate.
The two equivalent reformulations of \eqref{eq:original} differ only in the constraints \eqref{eq:diff1} and \eqref{eq:diff2}. 
The reformulation of \eqref{eq:original} in \eqref{eq:partial} allows us to solve \eqref{eq:original} by using Algorithm \ref{alg:nfeasible}. Algorithm \ref{alg:nfeasible} provides a distributed solution to the Economic Dispatch problem in \eqref{eq:original} if suppliers measure the used power $\{\mathbf{p}_i^{[k]}\},$ set the price $\{\boldsymbol{\lambda}^{[k]}\}$ and broadcast the new set of prices. 
This algorithm bears great resemblance to the OD3 algorithm proposed in \cite{tracking} with the exception of the box constraints and primal-dual inertia terms in \eqref{eq:inertia}.  However, in \cite{tracking}, the box constraints are omitted so as to obtain Q-linear convergence using Dual Decomposition.

 While Algorithm \ref{alg:nfeasible} or OD3 require a much smaller communication effort than Algorithm \ref{alg:feasible} (to be introduced), the box constraints are essential to the problem as they represent generation and consumption limits on generators and consumers respectively. Based on \cite{tracking} the amount by which \eqref{eq:powerbalance} is violated, i.e. $\mathbf{e}^{[k]} \triangleq \sum_{i=1}^N \mathbf{p}_i^{[k]} - \mathbf{P}^{[k]}$ can be bounded in norm. However, we can not know before-hand the sign of the components of $\mathbf{e}^{[k]},$ which denote an excess or shortage of power supply. This is critical since a widespread system breakdown (blackouts) can occur if the aggregate power consumption exceeds the supply capacity \cite{SindriDA}. Further, to the best of our knowledge, the convergence rate that ADMM achieves for Algorithm \ref{alg:nfeasible} does not suffice to claim that the iterates will remain at a bounded distance of the optimum point without the statement being trivial (i.e. due to compactness of the feasible set).
Under the assumptions given in Section 3, Algorithm 1 converges R-linearly in the static case. This can be proven using the results in \cite{luoADMM}.

The reformulation in \eqref{eq:total} allows us to solve \eqref{eq:original} by using Algorithm \ref{alg:feasible}. In this case, the problem each node solves, i.e. the iterate in \eqref{eq:p-iterate} consists of minimizing a strongly convex function with no constraints. This allows us to establish Q-linear convergence via \cite{linearADMM}. However, all constraints must be left for the iterate in \eqref{eq:coop}. This has the downside of requiring information exchange in order to be able to obtain the iterate \eqref{eq:coop}. The specifics of how this is done are included in Algorithm 3. For convenience, let us define $\mathcal{Q}^{[k]}_i$ as the set of vectors fulfilling constraint $i$ in \eqref{eq:diff2} for time $k.$ Further, let $[\mathbf{x}]_{\mathcal{Q}^{[k]}_i}$ denote the projection of $\mathbf{x}$ over $\mathcal{Q}^{[k]}_i.$ In particular, we require, for a bus network, the exchange of $RN + 2\sum_{j=1}^R |\mathcal{T}(j)|$ (cf. \eqref{eq:cf}) real quantities and the broadcasting of a binary vector of size $R.$
If we were to directly send the limits of the box constraints to solve the iterate \eqref{eq:coop} we would require the exchange of $4RN$ real quantities. Hence, in the \emph{worse case scenario} we will transmit $4RN$ real quantities and a binary vector of size $R.$ 
\begin{figure}
\vspace{-20pt}
\begin{algorithm}[H]
\caption{ADMM applied to \eqref{eq:partial}: Partial feasibility}
\label{alg:nfeasible}
\begin{algorithmic}[1]
	\State Initialize $\boldsymbol{\lambda}^{[0]} \,$ and set $k = 0.$ 
	\State Each node $i$ obtains $f_i(\cdot,k+1),$ $\underline{p}_i^{[k+1]}$ and $\overline{p}_i^{[k+1]}.$ The suppliers obtain $\mathbf{P}^{[k+1]}$
	\State Each user computes:
	\begin{subequations}
	\begin{align}
	 \mathbf{p}_i(k+1) := & \min_{\mathbf{p}_c} f_i(\mathbf{p}_i,k+1) + (\boldsymbol{\lambda}^{[k]} - \boldsymbol{\lambda}^{[k-1]})^T\mathbf{p}_i + 
	\label{eq:inertia} \\ & + \frac{\rho}{2}\|\mathbf{p}_i - \mathbf{p}_i^{[k]}\|^2 \nonumber \\
	 & \text{s.t.} \quad \underline{\mathbf{p}}_i^{[k+1]} \leq \mathbf{p}_i \leq \overline{\mathbf{p}}_i^{[k+1]}
	\end{align}
	\end{subequations}
	\State  The suppliers measure $\sum_{i=1}^N \mathbf{p}_i^{[k+1]} - \mathbf{P}^{[k+1]}$ and compute the price
	\begin{equation}
	\boldsymbol{\lambda}^{[k+1]} := \boldsymbol{\lambda}^{[k]} + \frac{\rho}{N} (\sum_{i=1}^N \mathbf{p}_i^{[k]} - \mathbf{P}^{[k+1]}).
	\end{equation}
	\State The suppliers broadcast the prices $\boldsymbol{\lambda}^{[k+1]}.$
	\State Set $k = k +1 $ 
\end{algorithmic}
\end{algorithm}
\vspace{-20pt}
\end{figure}

 Note that at no moment do we require the exchange of information regarding the objective functions. 
Algorithm \ref{alg:feasible} provides, by solving \eqref{eq:coop}, a feasible iterate at each iteration. The power balance constraint can be replaced by an inequality constraint if we do not require all the power to be used, i.e., we have storing devices. Note that, just as for Algorithm \ref{alg:nfeasible}, for the results in \cite{tracking} the direction of the constraint violation, i.e., the sign of each of the coordinates in $\mathbf{e}^{[k]},$ can not be determined and therefore the power balance constraint may not be fulfilled even if it is replaced by an inequality. Further, by replacing the power balancing constraint by an inequality constraint we can essentially turn the problem into a resource allocation problem with strongly concave utilities.
\begin{figure}
\vspace{-20pt}
\begin{algorithm}[H]
\caption{ADMM applied to \eqref{eq:total}: Total feasibility}
\label{alg:feasible}
\begin{algorithmic}[1]
	\State Initialize $\{\mathbf{p}_i^{[0]}\}_{i=1}^N$ and $\{\boldsymbol{\lambda}^{[0]}_i\}_{i=1}^N.$ Set $k = 0.$
	\State Each node $i$ obtains $f_i(\cdot,k+1),$ $\underline{p}_i^{[k+1]}$ and $\overline{p}_i^{[k+1]}.$ The system operator obtains $\mathbf{P}^{[k+1]}.$
	\State The system operator and the users cooperatively solve:
	\begin{subequations} \label{eq:coop}
	\begin{align}
	& \{\mathbf{q}_i^{[k+1]}\}_{i=1}^N =   \underset{\{\mathbf{q}_i\}_{i=1}^N}{\text{argmin}} \quad \frac{1}{2}\sum_{i=1}^N \|\mathbf{q}_i - (\mathbf{p}_i^{[k]} + \frac{\boldsymbol{\lambda}^{[k]}_i}{\rho})\|^2 \\
	& \text{s.t.} \qquad \underline{\mathbf{p}}_i^{[k + 1]} \leq \mathbf{q}_i \leq \overline{\mathbf{p}}_i^{[k+1]} \\
	& \quad \qquad \sum_{i=1}^N \mathbf{q}_i = \mathbf{P}^{[k+1]}
	\end{align}
	\end{subequations}
	\textbf{See Algorithm 3 for how to solve \eqref{eq:coop} cooperatively.}
	\State Each node computes:
	\begin{align}
	\mathbf{p}_i^{[k+1]} = \underset{\mathbf{p}_i}{\text{argmin}} \, f_i(\mathbf{p}_i,k+1) + \boldsymbol{\lambda}_i^{[k]T}\mathbf{p}_i  \label{eq:p-iterate}\\
	 \qquad + \frac{\rho}{2}\|\mathbf{p}_i - \mathbf{q}_i^{[k]}\|^2 \nonumber
	\end{align}
	\begin{equation}
	\boldsymbol{\lambda}_i^{[k+1]} = \boldsymbol{\lambda}_i^{[k]} + \rho(\mathbf{p}_i^{[k+1]}-\mathbf{q}_i^{[k+1]})
	\end{equation}
	\end{algorithmic}
\end{algorithm}
\vspace{-20pt}
\end{figure}
\section{Tracking Statement}
Consider the following assumptions.
\begin{assumption}[Uniform bounds on the curvature] \label{Assumption:curvature} The objective functions $f_i^{[k]}(\mathbf{p}_i) \triangleq C_i^{[k]}(\mathbf{p}_i)$ or $f_j^{[k]} (\mathbf{p}_j) \triangleq - U_j^{[k]}(\mathbf{p}_j)$ are $\sigma -$strongly convex and their gradients are $L -$Lipschitz continuous for all $k$, $i=1,\hdots,M$ and $j=M+1,\hdots,N$.
\end{assumption}

\begin{assumption}[Feasibility] \label{Assumption:feasibility} The optimization problem \eqref{eq:original} is feasible for all $k$. Further, it holds that $\sum_{i =1}^N \underline{\mathbf{p}}_i^{[k]} < \mathbf{P}^{[k]} < \sum_{i=1}^N \overline{\mathbf{p}}_i^{[k]}.$
\end{assumption}
\begin{assumption}[Bounded dynamics] \label{Assumption:dynamics} Let $\mathbf{p}^{\star[k]}$ denote the optimal point of \eqref{eq:original}. Then, the drift quantities are bounded, i.e. $\|\mathbf{p}^{\star[k]} - \mathbf{p}^{\star[k+1]}\| \leq \Delta {p}^{\star}$ and $\|\boldsymbol{\lambda}^{\star[k]} - \boldsymbol{\lambda}^{\star[k+1]}\| \leq \Delta {\lambda}^{\star},$
where the optimal dual multiplier  $\boldsymbol{\lambda}^{\star [k]}$ may correspond to either the constraint in the RHS of \eqref{eq:consistency_partial} or \eqref{eq:consistency} which become relevant when we reformulate the problem so as to solve it using ADMM. 
\end{assumption}
\begin{assumption}[Network connectivity] \label{Assumption:network} All nodes can reach the system operator, i.e. the network if fully connected.
\end{assumption}

Assumption \ref{Assumption:curvature} is a standard assumption for convex optimization methods that achieve Q-linear convergence rates; while linear convergence has been established under milder conditions \cite{strongconvex}, to the best of the authors' knowledge this still require appropriate step size selection to achieve such rates. On the contrary, ADMM achieves linear convergence rates (albeit the specific rate will vary with the step size) regardless of choice of step-size \cite{linearADMM}.

Assumption \ref{Assumption:feasibility} guarantees that not all generators and consumers will be pushed to the limit of their capabilities in order to fulfill the power balance constraint.
Assumption \ref{Assumption:feasibility} guarantees that the linear independence constraint qualification (LICQ) holds at the optimal point, implying uniqueness of the optimal dual multipliers \cite{unique}. In the static scenario uniqueness of the dual multipliers may not be a concern since it is sufficient to establish convergence to a KKT point. However, if the dual multipliers are allowed to move and they are not unique, further requirements are needed in order to define the drift of the multipliers. 

Assumption \ref{Assumption:dynamics} establishes a bound on the optimal power allocation and optimal multipliers from one iterate to the next. It can be shown that the quantities $\Delta \mathbf{p}^{\star}$ and $\Delta \boldsymbol{\lambda}^{\star}$ will remain bounded as long as the gradient drift is bounded \cite{trackADMM}, i.e., $ \exists \, \Delta f < \infty \text{ s.t. }\sum_i\|\nabla f_i(\mathbf{p}^{\star[k]},k) - \nabla f_i(\mathbf{p}^{\star[k+1]},k+1)\| \leq \Delta f.$

Assumption \ref{Assumption:network} is required so that the required information always reaches every node.
Before making the main statement we need to establish uniqueness of the optimal dual multipliers. This is done to guarantee that the multiplier drift $\|\boldsymbol{\lambda}^{\star [k]} - \boldsymbol{\lambda}^{\star[k+1]}\|$ is well defined.
\begin{lemma}
Under Assumptions 1-4 the optimal dual multipliers $\boldsymbol{\lambda}^{\star[k]}$ associated to \eqref{eq:consistency} are unique for each $k.$
\end{lemma}

\emph{Proof Sketch:}  All dual multipliers associated to the constraints in \eqref{eq:original} can be shown to be unique using the result in \cite{unique}. Further, by writing the optimality conditions of \eqref{eq:total} and \eqref{eq:original} we can show that $\boldsymbol{\lambda}^{\star[k]}$ is a linear combination of the optimal dual multipliers of \eqref{eq:original} implying the uniqueness of $\boldsymbol{\lambda}^{\star[k]}.$ \QEDB

Let $\mathbf{q}^{[k]},\,\mathbf{p}^{[k]},\,\boldsymbol{\lambda}^{[k]}$ be concatenations of $\{\mathbf{q}_i^{[k]}\}_{i=1}^N,\{\mathbf{p}_i^{[k]}\}_{i=1}^N$ and $\{\boldsymbol{\lambda}_i^{[k]}\}_{i=1}^N$ respectively.
\begin{theorem}
Under Assumptions \ref{Assumption:curvature}-\ref{Assumption:network} Algorithm \ref{alg:feasible} generates a sequence of iterates $\{\mathbf{q}^{[k]},\mathbf{p}^{[k]},\boldsymbol{\lambda}^{[k]}\}$ that fulfills 
\begin{equation}
\underset{k \to \infty}{ \text{lim sup}} \quad \|\mathbf{p}^{[k]} - \mathbf{p}^{\star[k]}\| \leq c_1, \text{and }
\end{equation}
\begin{equation}
\label{eq:usefulbound}
\underset{k \to \infty}{\text{lim sup}} \quad \|\mathbf{q}^{[k]}-\mathbf{q}^{\star{[k]}}\| \leq c_2,
\end{equation}
where $c_1 \triangleq \frac{g}{\sqrt{1+\delta_{\text{max}}}-1}, $ $c_2 \triangleq 3c_1^2 + \frac{1}{\rho}g^2 + \frac{3}{\sqrt{\rho}}c_1g ,$  $\delta_{\text{max}} \triangleq \frac{1}{\sqrt{L/\sigma}}$ and $g \triangleq \sqrt{\rho(\Delta \mathbf{p}^{\star})^2 + \frac{1}{\rho}(\Delta \boldsymbol{\lambda}^{\star})^2}.$ Note that the sequence $\{\mathbf{q}^{[k]}\}$ is always primal feasible.
\end{theorem}
\begin{figure}
\vspace{-20pt}
\begin{algorithm}[H]
\label{alg:projection}
\caption{Cooperative projection}
\begin{algorithmic}[1]
\State Each user computes $\mathbf{m}_i = [\mathbf{p}_i^{[k]} + \frac{\boldsymbol{\lambda}_i^{[k]}}{\rho}]_{\mathcal{Q}_i^{[k+1]}} - (\mathbf{p}_i^{[k+1]} + \frac{\boldsymbol{\lambda}_i^{[k]}}{\rho})$
\State Each user forwards information such that the system operator receives $\sum_{i=1}^N [\mathbf{p}_i^{[k]} + \frac{\boldsymbol{\lambda}_i^{[k]}}{\rho}]_{\mathcal{Q}_i^{[k+1]}}$
\State The system operator computes $\mathbf{d} =- \sum_{i=1}^N [\mathbf{p}_i^{[k]} + \frac{\boldsymbol{\lambda}_i^{[k]}}{\rho}]_{\mathcal{Q}_i^{[k+1]}} + \mathbf{P}^{[k+1]}$ and broadcasts $\text{sign}(\mathbf{d}).$
\State Each user then:
\If{$\text{sign}(\mathbf{d})(j) =\text{sign}( \mathbf{m}_i(j)) || \mathbf{m}_i(j)=0$}
\If{$\text{sign}(\mathbf{m}_i(j)) > 0$}
\State The user sends $(\mathbf{m}_i(j),\overline{x}_{ij})$ to the operator, where 
\begin{equation*}
\overline{x}_{ij} \triangleq \underline{\mathbf{p}}_i - [\mathbf{p}_i^{[k]} + \frac{\boldsymbol{\lambda}_i^{[k]}}{\rho}]_{\mathcal{Q}_i^{[k+1]}}(j).
\end{equation*}
\EndIf
\If{$\text{sign}(\mathbf{m}_i(j)) < 0$}
\State The user sends $(\mathbf{m}_i(j),\underline{x}_ij)$ to the operator, where $$\underline{x}_{ij} \triangleq \underline{\mathbf{p}}_i - [\mathbf{p}_i^{[k]} + \frac{\boldsymbol{\lambda}_i^{[k]}}{\rho}]_{\mathcal{Q}_i^{[k+1]}}(j).$$
\EndIf
\EndIf 
 \begin{equation}
\mathcal{T}(j) \triangleq \text{set of nodes transmitting regarding component }j. \label{eq:cf}
\end{equation}
\State The system operator solves 
\begin{subequations}
\begin{align}
\label{eq:coop_simple}
 \underset{\{\Delta q_{ij}\}_{i \in \mathcal{T}_2(j),j}}{\text{min}} & \sum_{ij}\| \Delta q_{ij} + \mathbf{m}_i(j)\|^2 \\
 \text{s.t.} \quad  & \underbrace{0 \leq \Delta q_{ij} \leq \overline{x}_{ij}\, i \in \mathcal{T}(j)}_{( \text{if } \text{sign}(\mathbf{d}(j)) > 0} \\
 & \underbrace{ \underline{x}_{ij} \leq \Delta q_{ij} \leq 0,\, i \in \mathcal{T}(j)}_{\text{if } \text{sign}(\mathbf{d}(j)) < 0} \\
  & \sum_{i \in \mathcal{T}(j)}\Delta q_{ij} = \mathbf{d}(j) 
\end{align}
\end{subequations}
\State The system operator sends to each node in $\mathcal{T}_(j)$ the required amount of movement, i.e. $\Delta\mathbf{q}_ij\,\,\forall i \in \mathcal{T}(j).$
\State Each node $j \in \mathcal{T}_2$ computes $\mathbf{q}_i^{[k+1]}(j) = [\mathbf{p}_i^{[k]}(j) + \frac{\boldsymbol{\lambda}_i^{[k]}(j)}{\rho}]_{\mathcal{Q}_i^{[k+1]}} + \Delta q_{ij}.$
\end{algorithmic}
\end{algorithm}
\vspace{-30pt}
\end{figure}

\begin{proof}
\vspace{-10pt}
This proof consists of two parts. Part 1 corresponds to the tracking statement based on ADMM's Q-linear convergence \cite{linearADMM}. The second part proves that Algorithm 3 actually solves the optimization problem in \eqref{eq:coop}.
\emph{Part 1:}
The first statement of the theorem follows from statements in \cite{linearADMM} and \cite{trackADMM}. In particular, the problem in \ref{eq:total} can be seen as an instance of \emph{Scenario 1} in \cite{linearADMM} for which Deng and Yin establish Q-linear convergence for $\{\mathbf{p}^{[k]},\boldsymbol{\lambda}^{[k]}\}$ if the problem is kept static. Let $\mathbf{u}^{[k]} \triangleq [\mathbf{p}^{[k]T},\boldsymbol{\lambda}^{[k]T}]^T,$ $\mathbf{u}^{\star} \triangleq [\mathbf{p}^{\star T},\boldsymbol{\lambda}^{\star T}]^T,$ $\mathbf{G} = \begin{pmatrix}  \rho \mathbf{I} & \mathbf{0} \\  \mathbf{0} & \frac{1}{\rho}\mathbf{I} \end{pmatrix}$ and $\|\cdot\|_{\mathbf{G}}$ be the semi-norm induced by $\mathbf{G,}$ from Theorem 3.1 and Corollary 3.1 in \cite{linearADMM} we have
\begin{equation}
\|\mathbf{u}^{[k+1]} - \mathbf{u}^{\star}\|_{\mathbf{G}}^2 \leq \frac{1}{1+\delta_{\text{max}}} \|\mathbf{u}^{[k]} - \mathbf{u}^{\star}\|_{\mathbf{G}}^2,
\end{equation}
with $\delta_{\text{max}} = \frac{1}{\sqrt{L / \sigma}}$ corresponding to selecting $\rho = \sqrt{\frac{L\sigma}{(NR)}}.$
For the dynamic case, this implies that
\begin{equation}
\|\mathbf{u}^{[k+1]} - \mathbf{u}^{\star [k+1]}\|_{\mathbf{G}}^2 \leq \frac{1}{1+\delta_{\text{max}}}\|\mathbf{u}^{[k]} - \mathbf{u}^{\star [k+1]}\|_{\mathbf{G}}^2,
\end{equation}
where $\mathbf{u}^{\star [k+1]}$ is now parametrized with an iteration number so as to indicate that the optimal primal-dual point moves over time. By taking square root, using the triangle inequality and evaluating the bound recursively, we obtain
\begin{align}
&\|\mathbf{u}^{[k+1]}-\mathbf{u}^{\star [k+1]}\|_{\mathbf{G}} \leq \left( \frac{1}{\sqrt{1+\delta_{\text{max}}}} \right)^{k+1}\|\mathbf{u}^{[0]}-\mathbf{u}^{\star [0]}\| \nonumber \\ + &\sum_{i=0}^{k}\left(\frac{1}{\sqrt{1+\delta_{\text{max}}}}\right)^{k-i+1}\|\mathbf{u}^{\star [i]} - \mathbf{u}^{\star [i+1]}\|_{\mathbf{G}},
\end{align}
where $\|\mathbf{u}^{\star[k+1]} - \mathbf{u}^{\star[k]}\|_{\mathbf{G}} \leq g.$
By taking the limits
\begin{equation}
\underset{k \to \infty}{\text{lim sup}} \quad \|\mathbf{u}^{[k+1]} - \mathbf{u}^{[k]}\|_{\mathbf{G}}  \leq \frac{g}{\sqrt{1+\delta_{\text{max}}} -1}.
\end{equation}
Up until here the procedure is but a simplified version of that in \cite{trackADMM}. We now proceed to derive the bound \eqref{eq:usefulbound} which is more interesting than the previous since it concerns the primal feasible iterate. Problem \eqref{eq:coop} finds the projection of $\mathbf{p}^{[k]} + \frac{\boldsymbol{\lambda}^{[k]}}{\rho}$ on the compact polyhedral set $\mathcal{P}^{[k+1]}$ defined by  $ \mathcal{P}^{[k+1]} \triangleq \{\mathbf{q} :\underline{\mathbf{p}}^{[k+1]} \leq \mathbf{q} \leq \overline{\mathbf{p}}^{[k+1]},\, \sum_{i=1}^N \mathbf{q}_i = \mathbf{P}^{[k+1]}\}.$ Hence we will write $\mathbf{q}_i^{[k+1]} = [\mathbf{p}_i^{[k]} + \frac{\boldsymbol{\lambda}_i}{\rho}]_{\mathcal{P}^{[k+1]}}.$ An intermediate step to showing that \eqref{eq:usefulbound} is true is that $\mathbf{q}_i^{\star [k+1]} = [\mathbf{p}_i^{\star [k+1]} + \frac{\boldsymbol{\lambda}^{\star [k+1]}}{\rho}]_{\mathcal{P}^{[k+1]}}$ for any $\rho > 0$ which can be done by writing the optimality conditions of \eqref{eq:original} and \eqref{eq:total}. We have already established that
\begin{align}
\label{eq:lambdaupper}
& \frac{1}{\rho}\|\boldsymbol{\lambda}^{[k+1]} - \boldsymbol{\lambda}^{\star[k+1]}\|^2  \leq \\ &\left( \frac{1}{\sqrt{1+\delta_{\text{max}}}} \|\mathbf{u}^{[k]} - \mathbf{u}^{\star[k]}\| + g \right)^2
- \rho \|\mathbf{p}^{[k+1]} - \mathbf{p}^{\star [k+1]}\|^2. \nonumber
\end{align}
Note that $\boldsymbol{\lambda}^{[k+1]} = \boldsymbol{\lambda}^{[k]} + \rho (\mathbf{p}^{[k+1]} - \mathbf{q}^{[k+1]})$ and $\boldsymbol{\lambda}^{\star [k+1]} = \boldsymbol{\lambda}^{\star [k+1]} + \rho(\mathbf{p}^{\star [k+1]} - \mathbf{q}^{\star [k+1]}),$ which means the LHS of \eqref{eq:lambdaupper} can be equivalently written as
\begin{align}
\rho \|\mathbf{q}^{[k+1]}-\mathbf{q}^{\star[k+1]}\|^2 +  \nonumber \\ \rho\left\Vert\mathbf{p}^{[k+1]}+ \frac{\boldsymbol{\lambda}^{[k]}}{\rho} - \left( \mathbf{p}^{\star [k+1]} + \frac{\boldsymbol{\lambda}^{\star [k+1]}}{\rho} \right)\right\Vert^2 + \nonumber\\
-2\rho\left( \mathbf{q}^{[k+1]}-\mathbf{q}^{\star[k+1]}\right)^T \\ \times \left(\mathbf{p}^{[k+1]}+ \frac{\boldsymbol{\lambda}^{[k]}}{\rho} - \left( \mathbf{p}^{\star [k+1]} + \frac{\boldsymbol{\lambda}^{\star [k+1]}}{\rho} \right) \right), \nonumber
\end{align}
where the last term can be lower bounded by $-2\rho \|\mathbf{p}^{[k+1]}+ \frac{\boldsymbol{\lambda}^{[k]}}{\rho} - \left( \mathbf{p}^{\star [k+1]} + \frac{\boldsymbol{\lambda}^{\star [k+1]}}{\rho} \right)\|^2$ because $\mathbf{q}^{[k+1]} = [\mathbf{p}^{[k+1]}+\frac{\boldsymbol{\lambda}^{[k]}}{\rho}]_{\mathcal{P}^{[k+1]}},$  $\mathbf{q}^{\star[k+1]} \triangleq [\mathbf{p}^{\star [k+1]} + \frac{\boldsymbol{\lambda}^{\star [k]}}{\rho}]_{\mathcal{P}^{[k+1]}}$ and $\mathcal{P}^{[k+1]}$ is convex. Through algebraic manipulation and using the triangular inequality  we obtain
\begin{align}
& \quad \|\mathbf{q}^{[k+1]} -\mathbf{q}^{\star [k+1]}\|^2 \leq \nonumber \\ & \quad \left( \frac{1}{\sqrt{1+\delta_{\text{max}}}}\left( \|\mathbf{u}^{[k]} - \mathbf{u}^{\star[k]}\| + g \right) \right)^2  \nonumber\\
& +\|\mathbf{u}^{[k]}-\mathbf{u}^{\star[k]}\|^2 + \frac{1}{\rho}(\Delta \boldsymbol{\lambda}^{\star})^2 + 
  \frac{2}{\sqrt{\rho}}\|\mathbf{u}^{[k]} - \mathbf{u}^{\star[k]}\|\Delta \boldsymbol{\lambda}^{\star}  \nonumber \\ & + \|\mathbf{u}^{[k+1]} - \mathbf{u}^{\star [k+1]}\|\|\mathbf{u}^{[k]}-\mathbf{u}^{[\star [k]]}\| \nonumber\\  & +  \quad \frac{\Delta \boldsymbol{\lambda}^{\star}}{\sqrt{\rho}}\|\mathbf{u}^{[k+1]} - \mathbf{u}^{\star[k+1]}\|.
\end{align}
Taking the limit on both sides yields 
\begin{align*}
\underset{k \to \infty}{\text{lim sup}} \quad \|\mathbf{q}^{[k+1]} - \mathbf{q}^{\star [k+1]}\|^2 \leq 3c_1^2 + \frac{1}{\rho}g^2 + \frac{3}{\sqrt{\rho}}c_1g 
\end{align*}
\emph{Part 2:} Algorithm 3 solves the optimization problem \eqref{eq:coop} by forcing the nodes to first obtain a feasible projection over their own box constraints. Then, the non-feasibility of the projections with regard to the power balance constraint \eqref{eq:powerbalance} is computed. It is then relevant to inspect in which direction the power balance constraint is violated. Depending on this direction, the variables that are in the lower or upper bound of their box constraints will be set to their optimal values already. Within the remaining variables, their value will be changed by the same amount, constraints permitting. We will now show that Algorithm 3 solves \eqref{eq:coop}. First of all, note that the solution to \eqref{eq:coop}, $\mathbf{q}^{[k+1]},$ can be written for each user as $\mathbf{q}^{[k+1]}_i(j) = [\mathbf{p}_i^{[k+1]} + \frac{\boldsymbol{\lambda}_i^{[k]}}{\rho}]_{\underline{\mathbf{p}}_i}^{\overline{\mathbf{p}}_i}(j) + \Delta q_{ij},$ where $\Delta q_{ij}$ indicates the deviation due to requiring the fulfillment of the power balance constraint. Problem \eqref{eq:coop} can then be equivalently written as
\begin{subequations}
\begin{align}
&\underset{\{\Delta q\}_{ij}}{\text{min}} & \sum_{i,j}\|\Delta q_{ij} + \mathbf{m}_{i}(j)\|^2 \\
& \text{s.t.} & \underline{x}_{ij} \leq \Delta q_{ij} \leq \overline{x}_{ij} \\
& & \sum_{i = 1}^N \Delta q_{ij} = \mathbf{d}(j). \label{eq:feas}
\end{align}
\end{subequations}
Whenever the quantity $\mathbf{d}(j) > 0,$ and $\overline{x}_{ij} = 0,$ we have that $\Delta q_{ij} \leq 0$ and $\mathbf{m}_i(j) \leq 0$ (cf. Step 1 in Algorithm 3). This implies that setting $\Delta q_{ij} \neq 0$ does not bring us closer to fulfilling \eqref{eq:feas} while it increases the value of the objective function. The analogous can be argued for $\mathbf{d}(j) < 0.$ Hence, solving \eqref{eq:coop_simple} is equivalent to solving \eqref{eq:coop}.
\end{proof}
\vspace{-12pt}
\section{Numerical Experiments}

\begin{figure}[H]
\vspace{-12pt}
\includegraphics[width = \columnwidth]{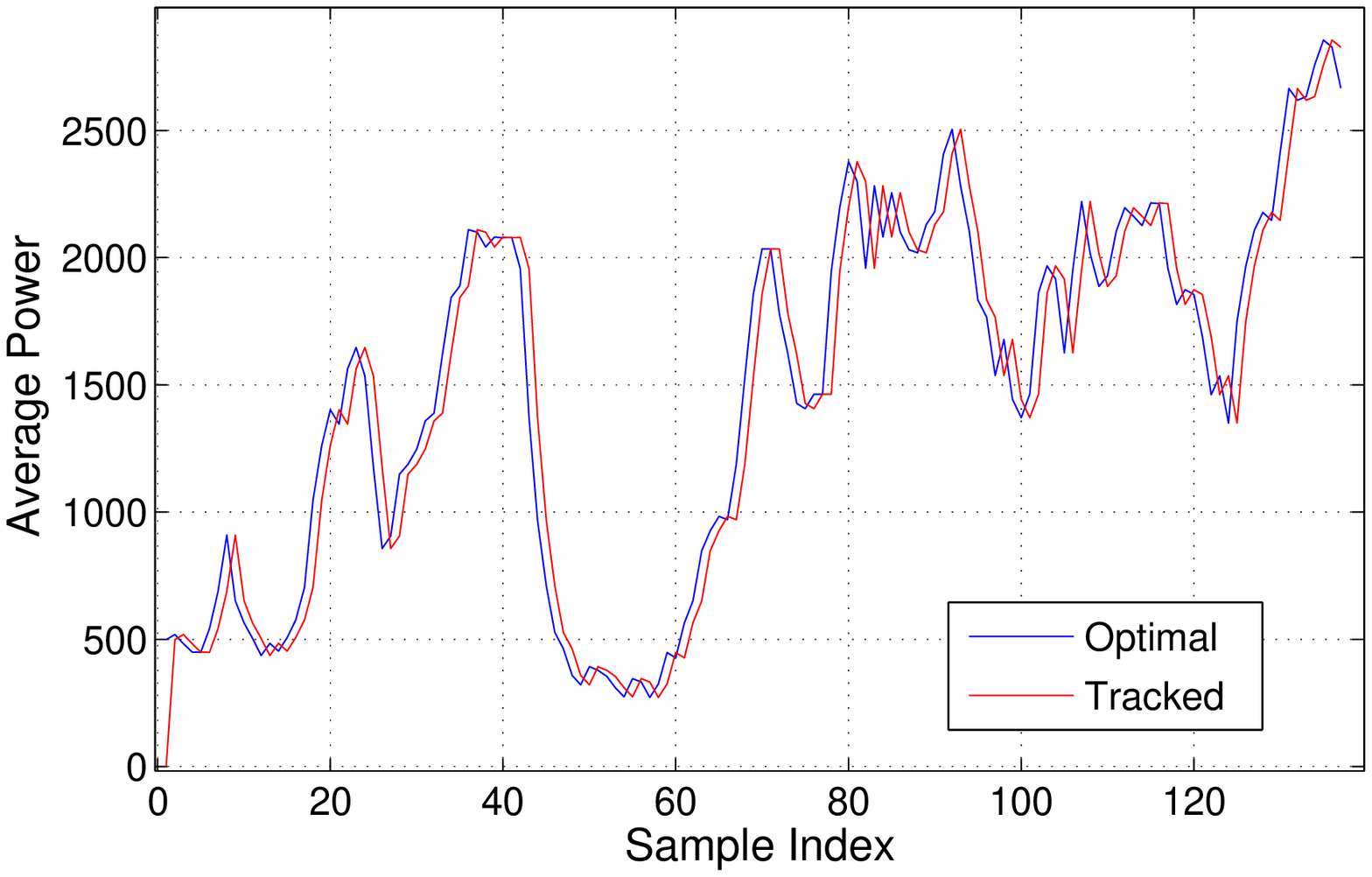}
\vspace{-10pt}
\caption{$\mathbf{p}^{\star[k]T}\mathbf{1}$(blue) and $\mathbf{q}^{[k]T}\mathbf{1}$ (red) generated by Algorithm \ref{alg:feasible}. Step size $\rho = 10.$}
\label{fig:only_figure}
\end{figure}

In this section we demonstrate numerically that the tracked solution remains close to the optimal solution when using Algorithm \ref{alg:feasible}. We use as in \cite{tracking} power generation data from the \emph{IEsystem operator Canada} Independent Electricity System Operator. We also consider a 10 user network with a single supplier, i.e. $N=10$ and $R=1.$ For each user $i=1,\hdots,10$, the cost function takes the form $(p_i - d_i^{[k]})^2$ where the demand $d_i^{[k]}$ is recursively updated as $d_i^{[k+1]} = [d_i^{[k]} + n^{[k]}]_+$ with $n^{[k]} \sim \mathcal{N}(0,1)$ and $d_i^{[0]} = 2.$ $P^{[k]}$ corresponds to the supply of aggregate power provided by renewable sources: biofuel, wind and solar, obtained in 5-minute intervals. The box constraints are set to $\underline{p}_i = 0$ and $\overline{p}_i = P^{[k]}$ for $i=1,\hdots,10.$ Fig. \ref{fig:only_figure} depicts and verifies that aggregate primal feasible iterates \eqref{eq:coop} track the optimal aggregate power. 

\section{Conclusions}
We have considered an economic dispatch problem where the utilities and constraints vary over time. As a difference to \cite{tracking} our analysis takes into account the problem's box constraints. Further, at the expense of increased information exchange, we are capable of providing a feasible solution at each iteration.
\bibliographystyle{IEEEbib}

\end{document}